\documentclass[12pt]{article}
\input cyracc.def

\usepackage{fullpage,amsmath,amscd,amsthm,amssymb,amsxtra,latexsym,epsfig}

\newcommand{\gp}{\mathfrak p}
\newtheorem{theorem}{Theorem}

\begin{document}

\centerline{\bf \Large An algorithm for computing the integral closure 
} 
\vskip10pt
\centerline{Theo de Jong}
\vskip20pt
In this paper we describe an algorithm for computing
the normalization for certain rings. This quite general
algorithm is essentially due to Grauert and Remmert
\cite{gr1}, \cite{gr2}, 
and seems to have escaped the attention
of the computer algebra specialists until now. 
As I am not a computer algebra specialist myself, I do not
know whether this algorithm is fast. 
Grauert and Remmert proved a normality
criterium in order to give a simple and nice proof of a theorem of
Oka, who proved that the non-normal points of 
an analytic space is an analytic space itself.
We just reformulate their results,  so that it 
better suits our purposes.\\

To fix the notation, let $R$ be a Noetherian 
and reduced ring, 
$\widetilde{R}$ the integral closure (also called normalization)
of $R$.
We consider the set:
\[
NNL:= \{ {\gp} \in {\rm Spec}(R) : R_{\gp} {\rm \; is\;  not \;normal} \}
\]
Here $NNL$ stands for non-normal locus. 
Let $I$ be an ideal of $R$
containing a nonzerodivisor. We have canonical inclusions:
\[
R \subset {\rm Hom}_R(I,I) \subset \widetilde{R}
\]
The first inclusion is the map which sends an element of
$R$ to multiplication with this element. The second inclusion
is  sending $\phi \in {\rm Hom}_R(I,I)$ to $\frac{\phi (f)}{f}$ for
any element $f \in I$ which is a nonzerodivisor of $R$. 
It is easily checked that the map is independent of
the choice of $f$.
That we in fact land in $\widetilde{R}$ can be found
in any textbook which has included integral closure as a topic. 

\begin{theorem} \cite{gr1} pp. 220-221, \cite{gr2}, pp. 125-127.
Assume that the ideal 
$I$ contains a nonzerodivisor, and  has the following property:
\[
NNL \subset V(I)
\]
where $V(I)= \{ {\gp} \in {\rm Spec} (R): I \subset {\gp} \}
$ denotes, as usual, the zero set of
$I$.
Suppose moreover that $I$  has the property 
\[
 {\rm Hom}_R(I,I) = {\rm Hom}_R(I,R) \cap \widetilde{R} \hspace{5mm} (*)
\]
Then one has the following normality criterium:
\[
R = {\rm Hom}_R(I,I) \iff R {\rm \; is \; normal}
\]
\end{theorem}

\begin{proof}
The implication $\Longleftarrow $ is trivial. To prove the
converse, let 
 $h = \frac{f}{g} \in \widetilde{R}$. 
Consider the following ideal in $R$
\[
 \{ g \in R: hg \in R \} 
\]
Its zero set is called the "pole set" of $h$:
\[
P(h):= \{ {\gp} \in {\rm Spec} (R): h \notin R_{\gp} \}
\]
It is immediate that $P(h) \subset NNL$. 
Let $J$ be the ideal of $P(h)$. There exists a $c > 0$ 
such that $h J^c \subset R$, by the Nullstellensatz.
By the 
Nullstellensatz again  $\sqrt{I} \subset J$.
Therefore there exists a $d > 0$ such that
$hI^d \subset R$. Let $d$ be minimal with this
property. We claim $d = 1$. Suppose the converse, i.e.
$d >1$. Then there exists  an $a \in I^{d-1}$ with
$ha \notin R$. Furthermore $ha \in \widetilde{R} $ 
and $(ha) I \subset R$.  
By assumption:
\[
R= \{ h \in \widetilde{R}: hI \subset I \}
= \{ h \in \widetilde{R}: hI \subset I \}
\]
so that $(ha) \in I$. Therefore $ha \in R$ after all,
a contradiction. 
\end{proof}
We have to find an ideal which satisfies condition $(*)$. 
This is provided by:

\begin{theorem}
Every radical ideal $I$ containing a nonzerodivisor safisfies
condition $(*)$.
\end{theorem}
\begin{proof}
The proof is in \cite{gr1} and \cite{gr2}, but because 
it is so nice and simple we give it here. Let 
$h \in \widetilde{R}$, so we have an equation:
\[
h^n = a_0 + a_1h + \ldots + a_{n-1}h^{n-1}; \hspace{5mm} a_i \in R
\]
If $hI \subset R$, then we have for all $f \in I$:
\[
(hf)^n = a_0f^n + a_1h + \ldots + a_{n-1}h^{n-1}f \in I
\]
As $I$ is supposed to be reduced it follows that $hf \in I$,
and that is what we had to prove. 
\end{proof}
These results give rise to the following algorithm:
\vskip10pt
\centerline{\bf \large ALGORITHM}
\vskip5pt
\noindent INPUT: A reduced noetherian ring $R$.\\
OUTPUT: The normalization $\widetilde{R}$ of $ R$.\\
STEP 1: Determine an ideal $I$ with $NNL \subset V(I)$.\\
STEP 2: Compute the radical $\sqrt{I}$ of $I$. Put
 $I := \sqrt{I}$.\\
STEP 3: Take an  $f \in I$, and compute $J := {\rm Ann}(f)$.
If $J=0$, GOTO STEP 5.\\
STEP 4: Put $R: = R/(f) \oplus R/J$ and GOTO STEP 1.\\
STEP 5: Compute ${\rm Hom}_R(I,I)$. If 
$R = {\rm Hom}_R(I,I)$ then put $\widetilde{R} := R $ and STOP.\\
STEP 6: Set $R:= {\rm Hom}_R(I,I)$ and GOTO STEP 1.\\

This algorithm stops exactly when the normalization
$\widetilde{R}$ is finitely generated as an $R$-module,
so for example for affine rings, due to a classical
result of E. Noether. \\

Some remarks are in order.

\begin{enumerate}
\item To determine an ideal $I$ with 
$NNL \subset V(I)$ one can take any $I$ which
contains the non-regular locus of $R$. 
\item Algorithms for computing the radical of an ideal are
described in \cite{ei}
and probably will take the longest in this algorithm.
\item In step 4, $R \subset R/(f) \oplus R/J$ is indeed
an inclusion, as $(f) \cap J =0$ follows from the assumption
$R$ is reduced.
This extension is also finite. To see this, take
a $g \in J$ with $f+g$ a nonzerodivisor (prime avoidance).
The extension:
\[
R \subset R/(f) \oplus R/(g) \cong R[X]/(X^2 -X, X(f+g)-f)
\]
obviously is finite. As $R \subset R/(f) \oplus R/J$
is a quotient of $ R/(f) \oplus R/(g)$ (in fact it is equal)
it is a finite extension too.

\item I think it should be possible to extend this
algorithm to an algorithm which computes the primary
decomposition for a radical ideal.  
\end{enumerate}

We finish with describing the ring structure 
of ${\rm Hom}_R(I,I)$, 
which essentially  is due to Catanese, see \cite{Ca} . Take generators
$u_0 := 1, u_1, \ldots , u_t$  of ${\rm Hom}_R(I,I)$ as $R$-module. 
Consider the map:
\[
 R\cdot X_0 \oplus  R\cdot X_1 \oplus \ldots 
\oplus R \cdot X_t \stackrel{\phi}{\longrightarrow} {\rm Hom}_R(I,I)
\]
\[
X_i \mapsto u_i
\]
Computing  the kernel of the map $\phi$ gives 
"linear equations":
\[
L_i = \sum_{j=0}^t{\alpha}_{ij}X_j = 0 \hspace{5mm} {\alpha}_{ij} \in 
R; \; i = 1, \ldots ,s
\]
Because ${\rm Hom}_R(I,I)$ is a ring, we have
that $u_iu_j$ for all $ 1 \leq i \leq j \leq t $ is in 
${\rm Hom}_R(I,I)$ again. We therefore can find elements
${\beta}_{ijk} \in R$ such that:
\[
u_iu_j = \sum_{k=0}^t{\beta}_{ijk}u_k
\]
giving us $\frac{t(t+1)}{2}$ "quadratic equations":
\[
Q_{ij} := X_iX_j - \sum_{k=0}^t{\beta}_{ijk}X_k
\]
For the easy proof of the following theorem we 
refer to Catanese \cite{Ca} .

\begin{theorem}
Put $X_0 = 1$, and consider
the ideal $J \subset R[X_1 , \ldots , X_t]$ generated
by the $L_i , i = 1, \ldots , s$ and  the 
$Q_{ij}$ for $ 1 \leq i \leq j \leq t$. 
Then we  have a ring isomorphism:
\[
{\rm Hom}_R(I,I) \cong R[X_1, \ldots , X_t]/J
\]
\end{theorem}
\vskip10pt
{\bf Acknowledgment} I am very grateful to Wolfram Decker for 
bringing to my attention that this algorithm was unknown
to the computer algebra specialists.

\end{document}